\documentclass[a4paper]{amsart}

\usepackage{amsmath,amsthm,amsfonts,amssymb}
\usepackage{graphicx}
\usepackage{latexsym}
\usepackage[all]{xy}
\usepackage{subfigure}
\usepackage[usenames,dvips]{color}  % for use of 68 colours

\usepackage{psfrag}  % allows to scan text in eps-files, to convert them to latex

\xyoption{matrix}
\xyoption{arrow}

\newtheorem{lemma}{Lemma}[section]

\newtheorem{theorem}[lemma]{Theorem}

\theoremstyle{definition}
\newtheorem*{remark}{Remark}
\newtheorem*{definition}{Definition}
\newtheorem{example}[lemma]{Example}

\newcommand{\QQ}{\mathbb Q}
\newcommand{\ZZ}{\mathbb Z}
\newcommand{\DD}{\mathcal D}
\newcommand{\CC}{\mathcal C}

\newcommand{\Hom}{\operatorname{Hom}\nolimits}

\begin{document}

\title[Cluster categories and diagonals]
{Cluster categories, $m$-cluster categories and diagonals in polygons}

\author[Baur]{Karin Baur}
\address{Department of Mathematics \\
ETH Z\"urich
R\"amistrasse 101 \\
8092 Z\"urich \\
Switzerland
}
\email{baur@math.ethz.ch}

\subjclass[2000]{Primary:   Secondary: }

\begin{abstract}
The goals of this expository article are on one hand 
to describe how to construct ($m$-)
cluster categories from triangulations (resp. from $m+2$-angulations) of polygons. 
On the other 
hand, we explain how to use translation quivers and their powers 
to obtain 
the $m$-cluster categories directly from the diagonals of 
a polygon. 
\end{abstract}

\maketitle

\tableofcontents

\section*{Introduction}

This expository article is the expanded version of a talk given at the 
conference at the 
Grenoble summer school "Geometric methods in representation theory" 
in July 2008. The goal of this talk was to explain how cluster categories 
and $m$-cluster categories can be described via diagonals and 
so-called $m$-diagonals in a polygon. 
And then how the latter can actually be described using the power of 
a translation quiver. 
The first section gives a very brief introduction to the theory of cluster algebras 
and cluster categories. It also introduces the notations used in the article. 
In Section 2, we explain the notions of a quiver 
given by the diagonals in a polygon and of the one given by $m$-diagonals. 
The results in this section are mainly due to Caldero-Chapoton-Schiffler (\cite{ccs}), 
to Schiffler (\cite{schi}) 
and to Baur-Marsh (\cite{bm},~\cite{bm2}). 
In the last section, we introduce the concept of the power of a translation 
quiver. Here, the results are from~\cite{bm},~\cite{bm2} and from 
the masters thesis of C. Ducrest (\cite{du08}).

\section{Cluster algebras and cluster categories}

Cluster algebras were introduced by Fomin and Zelevinsky (\cite{fz}) 
in order to provide 
an algebraic framework for the phenomena of total positivity and for 
the canonical bases of the quantized universal enveloping algebras.

We briefly illustrate the notion of total positivity: 
An $n\times n$ matrix is called {\em totally positive} if all its minors are 
positive. Originally, 
the term was used to describe matrices with non-negative minors: 
these matrices are nowadays called {\em totally non-negative}. 
In the 1930s, Gantmacher-Krein and I. Schoenberg have independently 
started investigating such matrices. One of the motivations was to 
estimate the number of real zeroes of a polynomial. 

Gantmacher showed that totally non-negative matrices 
have different real eigenvalues. 
The interest in total positivity was renewed in the 1990s when 
G. Lusztig extended the notion 
to reductive algebraic groups, cf.~\cite{lu}. 

\begin{example}
To illustrate the notion on a (non-reductive) example, let us consider the 
group of $3$ by $3$-matrices with $1$'s on the diagonal and zeroes below. 
If $U=\begin{small}
\begin{pmatrix} 1 & a & b \\ 0 & 1 & c \\ 0 & 0 & 1\end{pmatrix}\end{small}$ 
is such a matrix then $U$ is totally positive if $a,b,c>0$ and $ac-b>0$. 
One can check that it is actually enough to require $a>0,\ ac-b>0$ and $b>0$: 
the condition $c>0$ will follow automatically. 
Equivalently, the condition $a>0$ can be dropped  
and is automatically satisfied by the remaining conditions. 
So there is only a certain number of minors that need to be checked. 

Furthermore, if from the set $\{a>0,ac-b>0, b>0\}$ of conditions the first 
is omitted then we may replace it exactly with one other condition, namely 
with $c>0$, to obtain total positivity of the matrix. 
\end{example}

More generally, one shows that the minimal sets of minors to check all 
have the same cardinality. And it is often the case that if you remove one 
minor from such a minimal set, there exists exactly one other minor to replace 
it with. 

%%%%%%%%%%%%%%
%
\subsection{Cluster algebras} 
%
%%%%%%%%%%%%%%

A cluster algebra $\mathcal A\subset\QQ(u_1,\dots,u_n)$ of rank $n$ is 
an algebra with possibly infinitely many generators. 
These generators are called cluster variables; they are 
arranged in overlaping sets of the same cardinality $n$, the {\em clusters}. 
There are relations between the cluster variables, encoded in an $n\times n$ 
matrix, the {\em mutation matrix}. 
Through mutation, one element of a cluster is exchanged by exactly 
one other element and this exchange process is prescribed by the 
exchange matrix. 

If there are only finitely many generators, the cluster algebra is of finite 
type. Finite type cluster algebras have been classified by Fomin and 
Zelevinsky (\cite{fz-finite}). 
Their classification describes the finite type cluster algebras in terms 
of Dynkin diagrams. 

More concretely: a seed is a pair 
$(\underline{x},M)$ where $\underline{x}=\{x_1,\dots, x_n\}$ 
is a basis of $\QQ(u_1,\dots,u_n)$ and $M=(M_{ij})_{ij}$ is a sign skew symmetric 
$n\times n$-matrix with integer entries, called the exchange matrix. 
That means that the sign of $M_{ij}$ is the opposite of 
the sign of $M_{ji}$. 

Then one defines an involutive map $\mu_k$ (for $k\in\{1,\dots,n\})$ on the seeds, 
called the mutation in direction of $k$, through 
$\mu_k(\underline{x})=(x_1,\dots,\widehat{x}_k,x_k',\dots,x_n)$ where 
$x_k'$ is given by the relation 
$$
x_k\cdot x_k'=\prod_{
         \begin{array}{l}x_i\in\underline{x} \\ M_{ik}>0\end{array}
         } 
         x_i^{M_{ik}} + 
\prod_{
         \begin{array}{l}x_i\in\underline{x} \\ M_{ik}<0
        \end{array}
         } 
         x_i^{-M_{ik}}
$$
In a similar way, one defines $M'$ by 
$$
M_{ij}':=\left\{ 
              \begin{array}{ll}
              -M_{ij} & \mbox{if $i=k$ or $j=k$} \\
              M_{ij} + \frac{1}{2}(|M_{ik}|M_{kj} + M_{ik}|M_{kj}|) & \mbox{otherwise.}
              \end{array}
            \right.   
$$
and thus obtains 
$(\underline{x'},M')$ as $\mu_k((\underline{x},M))$ (the matrix $M'$ is also 
a sign skew symmetric $n\times n$-matrix 
over $\ZZ$). 
For more details we refer to Section 1 of the survey article 
\cite{bum} of Buan-Marsh. 
The $x_i$ obtained through successive mutations are the 
so-called cluster variables. 
The cluster algebra $\mathcal A={\mathcal A}(\underline{x},M)$ 
is then defined as the algebra generated by the cluster variables. 
There can be infinitely many of them. 
Fomin-Zelevinsky have shown in~\cite{fz-4} that $\mathcal A$ lies 
in $\ZZ[x_1^{\pm 1},\dots, x_n^{\pm 1}]$ (Laurent-phenomenon). 
First 
examples of cluster algebras are coordinate rings of SL$_2$, SL$_3$. 

The field of cluster algebras is a young and very dynamic field. Since 
its first introduction, there have been many different directions in its 
development. We only mention a few connections to other areas (in parentheses: 
the objects corresponding to the cluster variables): 
the theory of Teichm\"uller spaces (Penner coordinates), see work of 
Fock-Goncharov,~\cite{fg} and~\cite{fg2}; the representation 
theory of finite dimensional algebras (tilting modules), cf.~\cite{bmrrt}), 
triangulations of surfaces (diagonals), see the work~\cite{fst} of 
Fomin-Shapiro-Thurston. 

%%%%%%%%%%%%%%
%
\subsection{Cluster categories}
%
%%%%%%%%%%%%%%

Cluster categories were introduced independently in the work~\cite{bmrrt} 
of Buan-Marsh-Reineke-Reiten-Todorov,
and by Caldero-Chapoton-Schiffler,~\cite{ccs} 
to provide a categorification of the theory of cluster algebras. 

We will use the approach of~\cite{bmrrt} to describe 
cluster categories and will consider the approach of~\cite{ccs}Êlater, 
cf. Section~\ref{s:poly-diag}. 

Let $Q$ be a simply-laced Dynkin quiver (i.e. a quiver whose underlying 
graph is of type A, D or E). Let $k$ be an algebraically closed field and 
$kQ$ the path algebra of $Q$ (for more details, we refer to the 
lecture notes of M. Brion, \cite{br} in the same volume). 
Take the bounded derived category $\DD^b(kQ)$ of finitely generated 
$kQ$-modules (for details on $\DD^b(kQ)$ we refer to~\cite{ha88}) 
with shift denoted by $[1]$ and 
Auslander-Reiten translate denoted by $\tau$. By Happel (\cite{ha88}), 
the category $\DD^b(kQ)$ is triangulated, Krull-Schmidt, and 
has almost split sequences. 
To understand the category $\DD^b(kQ)$ 
it is helpful to study its Auslander-Reiten quiver: 
The Auslander-Reiten quiver of a category is a combinatorial tool which helps 
understanding the category. Its vertices are by definition the indecomposable modules 
up to isomorphism and the number of arrows between two points 
are given by the dimension of the space of irreducible maps between 
two representatives of the corresponding modules.

We now associate a quiver $\ZZ Q$ to $Q$. Its vertices are 
$(n,i)$ for $n\in\ZZ$, and where $i$ a vertex of $Q$. For every arrow $i\to j$ 
in $Q$ there are arrows $(n,i)\to (n,j)$ and $(n,j)\to(n+1,i)$ in $\ZZ Q$. 
So $\ZZ Q$ has the shape of a $\ZZ$-strip of copies of $Q$. 
Together with the map $\tau:(n,i)\to (n-1,i)$ ($n\in\ZZ$, $i$ a vertex of $Q$), 
$\ZZ Q$ is a stable translation 
quiver as defined by Riedtmann (see~\cite{rie}). For a precise definition of 
stable translation quivers we refer the reader to Section~\ref{s:powers} below.  
We illustrate $\ZZ Q$ in Example~\ref{ex:AR-A3}. 
Happel has shown in~\cite{ha88}, that the 
Auslander-Reiten quiver AR$(\DD^b(kQ))$ of $\DD^b(kQ)$ is just 
$\ZZ Q$. 
In particular, the category $\DD^b(kQ)$ is independent of the orientation of $Q$.

\begin{example}\label{ex:AR-A3}
Let $Q$ be a quiver of type A$_3$, 
$$
Q: \quad 
\includegraphics[scale=.4]{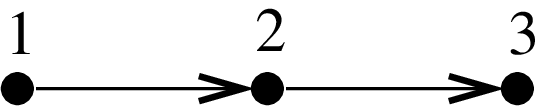}
$$
Then, $\ZZ Q$ has the shape 
$$
\includegraphics[scale=.5]{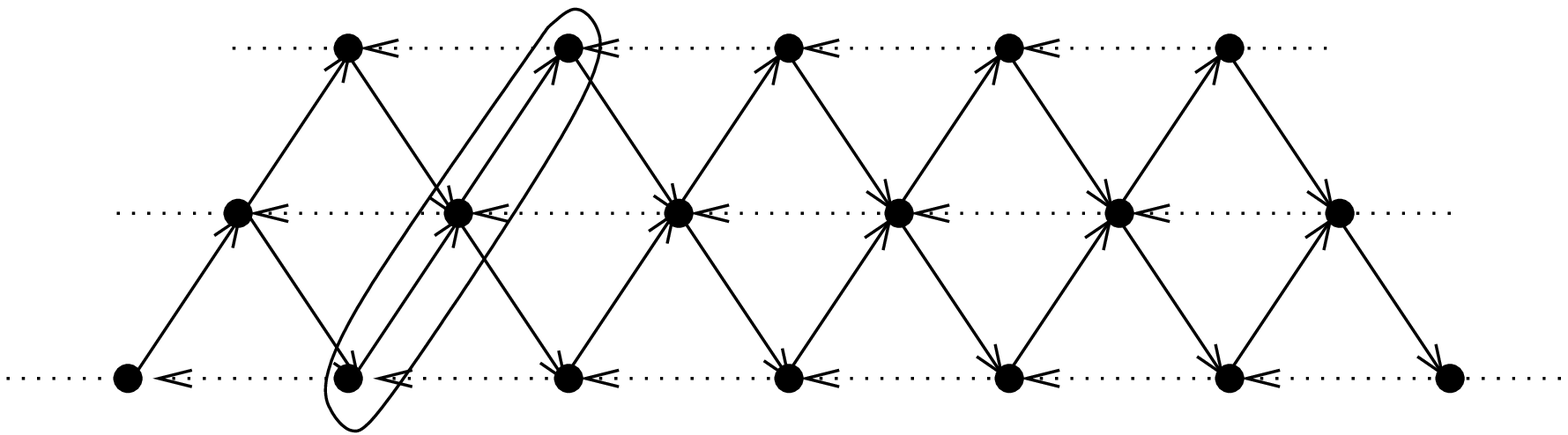}
$$
with one copy of the quiver $Q$ highlighted to show 
how it appears inside $\ZZ Q$. 
The dotted arrows indicate the Auslander-Reiten translate 
$\tau$ which sends each vertex to its leftmost neighbor. 
It is an auto-equivalence of the Auslander-Reiten quiver. 

On the other hand, the Auslander-Reiten quiver of the 
module category $kQ$-mod of finitely generated $kQ$-modules 
looks like a triangle: 
$$
\mbox{AR($kQ$-mod):}
\quad\quad 
\includegraphics[scale=.4]{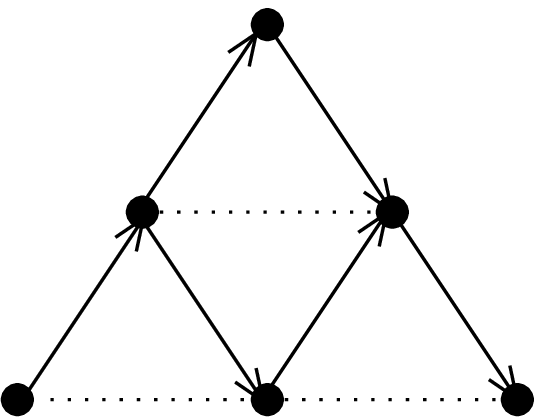}
$$

Observe that the infinite quiver $\ZZ Q$ can be viewed as being covered by 
copies of the Auslander-Reiten quiver of the module category $kQ$-mod, 
with additional arrows and dotted arrows introduced to connect the copies 
of the triangle of $kQ$-mod. 
With this picture in mind, we can describe the shift $[1]$ on 
AR$(\DD^b(kQ))$: it 
sends each vertex to the ``same'' vertex 
in the next copy of the triangle AR($kQ$-mod) to the right. 
\end{example}

Back to the general situation, where $Q$ is of simply-laced Dynkin type. 
The shift $[1]$ then is the auto-equivalence of 
AR$(\DD^b(kQ))$ which sends a vertex to the corresponding vertex 
in the next copy of the Auslander-Reiten quiver of the module category 
$kQ$-mod and the translation $\tau$ sends a vertex to its leftmost neighbor. 
As an abbreviation, we 
write $F$ for the auto-equivalence $\tau^{-1}\circ[1]$ of $\DD^b(kQ)$. 
Now we are ready to define the cluster category associated to 
$Q$. 

\begin{definition}
The {\em cluster category} $\CC:=\CC_Q:=\DD^b(kQ)/F$ {\em of type $Q$} is 
the orbit category whose objects are the $F$-orbits of objects of $\DD^b(kQ)$ 
and whose morphisms are given as follows: 
$$
\Hom_{\CC}(\widetilde{X},\widetilde{Y})=\bigsqcup_{i\in \ZZ} \Hom_{\DD^b(kQ)}(F^iX,Y)
$$
where $\widetilde{X}$ and $\widetilde{Y}$ are representatives of the 
$F$-orbits through $X$ and $Y$ respectively.  
\end{definition}

Note that for any pair of objects $X$, $Y$ of $\DD^b(kQ)$ 
there are only finitely many $i$ such that 
$\Hom_{\DD^b(kQ)}(F^iX,Y)$ is non-zero. 
The cluster category is Krull-Schmidt (\cite{bmrrt}), triangulated and 
Calabi-Yau of dimension 2 (\cite{ke}). 

The connection between $\CC_Q$ 
and the cluster algebra of the same type is given by the following result. 

\begin{theorem}[\cite{bmrrt}]
There is a bijection between 
the cluster variables of the cluster algebra of type A$_n$ (resp. 
D$_n$, E$_n$) and 
indecomposable objects of $\CC_Q$ where $Q$ is of type A$_{n}$ 
(resp. D$_n$, E$_n$). 
\end{theorem}

\begin{figure}
\begin{center}
  \psfragscanon
   \psfrag{X}{$X$}
   \psfrag{Y}{$Y$}
   \psfrag{Z}{$Z$}
\includegraphics[scale=.5]{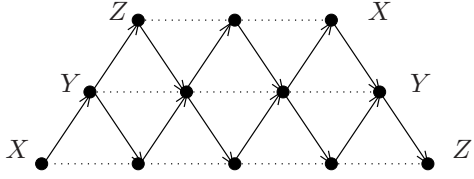}
\end{center}
\caption{The Auslander-Reiten quiver of $\CC$ for A$_3$}\label{fig:A3-AR}
\end{figure}

\begin{figure}
\begin{center}
  \psfragscanon
   \psfrag{X}{$X$}
   \psfrag{Y}{$Y$}
   \psfrag{U}{$U$}
   \psfrag{V}{$V$}
\includegraphics[scale=.5]{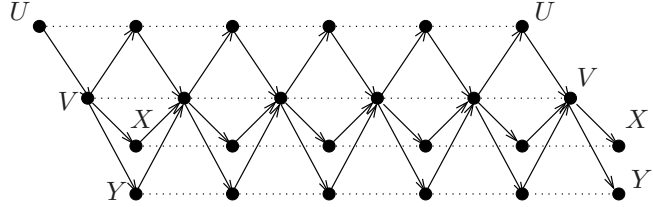}
\end{center}
\caption{The Auslander-Reiten quiver of $\CC$ for D$_4$}\label{fig:D4}
\end{figure}

To understand the cluster categories better, we consider its 
Auslander-Reiten quiver. By definition, it has the form of one copy 
of the module category, together with a copy of the quiver $Q$ 
(with additional arrows, dotted arrows), 
as illustrated in types A and D below (Figures~\ref{fig:A3-AR} and~\ref{fig:D4}). In 
particular, it is a finite quiver. 
In the pictures of Figures~\ref{fig:A3-AR} and~\ref{fig:D4}  
we have repeated one copy of the quiver $Q$ to indicate how the quivers are 
glued together: both quivers wrap around. 
The Auslander-Reiten quiver 
of type A$_n$ can be viewed as lying on a M\"obius strip and 
the one of type D$_n$ as lying on a cylinder. 

Let $G$ be the underlying graph of $Q$, $G$ of Dynkin type 
A, D or E. 
We recall that a famous result of P. Gabriel establishes a bijection between the indecomposable objects of $kQ$-mod (up to isomorphism) 
and the positive roots of the Lie algebra of type $G$. For a recent 
description of this result we refer to Section 5 of the lecture notes~\cite{kr} of 
H. Krause. 

There exists an analogous result for cluster categories. To state it we need 
to enlarge the set of roots considered: we add the negatives of the simple 
roots. 
An {\em almost positive root} of a Lie algebra of type $G$ is a positive root 
or the negative of a simple root. 

Buan et al. have shown in~\cite{bmrrt} that 
there is a bijection between the indecomposable 
objects of the cluster category $\CC_Q$ and the almost positive roots 
of the Lie algebra of type $Q$. 

%%%%%%%%%%%%%%%%%%%%
%
\subsection{The $m$-cluster category}
%
%%%%%%%%%%%%%%%%%%%%%

In 2005, Keller (\cite{ke}) has introduced the $m$-cluster categories 
as a natural generalisation of the cluster categories. 
Again, let $Q$ be a quiver whose underlying 
graph is of Dynkin type A, D or E. Let $[1]$ be the shift and $\tau$ the 
Auslander-Reiten translate as before. Let $F_m$ be the auto-equivalence 
$\tau^{-1}\circ[m]$ of $\DD^b(kQ)$, for $m\ge 1$. 

\begin{definition}
The {\em $m$-cluster category $\CC^m:=\CC_Q^m:=\DD^b(kQ)/F_m$ (of type $Q$)} 
is the orbit category with objects the $F_m$-orbits of objects of $\DD^b(kQ)$ and 
with morphisms 
$\Hom_{\CC^m}(\widetilde{X},\widetilde{Y})=\bigsqcup_{i\in \ZZ} 
\Hom_{\DD^b(kQ)}(F_m^iX,Y)$ where $\widetilde{X}$ and $\widetilde{Y}$ 
are representatives of the $F_m$-orbits of $X$ and $Y$ respectively. 
\end{definition}

\begin{figure}[h]
\begin{center}
  \psfragscanon
   \psfrag{X}{$X$}
   \psfrag{Y}{$Y$}
\includegraphics[scale=.5]{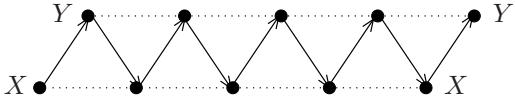}
\end{center}
\caption{The Auslander-Reiten quiver of $\CC^2$ for A$_2$}\label{fig:mA}
\end{figure}

Note that the Auslander-Reiten quiver of $\CC^m$ thus consists of 
$m$ copies of the Auslander-Reiten quiver 
of the module category $kQ$-mod and additionally, of a copy of the quiver $Q$, 
connected with additional (dotted) arrows. 
In the case of A$_n$, we observe that if 
$m$ is odd, the Auslander-Reiten quiver of $\CC_Q^m$ 
lies on a M\"obius strip, whereas if it $m$ is even, it lies on a cylinder. 
As an example, Figure~\ref{fig:mA} shows the Auslander-Reiten quiver 
of the $2$-cluster category of type A$_2$. Again, we have repeated a slice 
of the quiver to indicate how it wraps around. 

\begin{remark}
The $m$-cluster categories have very nice properties, analogously to the
properties of the cluster categories: 
$\CC^m$ is Krull-Schmidt (\cite{bmrrt}), triangulated and 
Calabi-Yau of dimension $m+1$ (\cite{ke}).  
\end{remark}

%%%%%%%%%%%%%%%%%%%%
%
\section{Polygons and diagonals}\label{s:poly-diag}
%
%%%%%%%%%%%%%%%%%%%%%

In the first part of this section we present the approach of 
Caldero, Chapoton and Schiffler,~\cite{ccs}, 
who described the cluster category of type A$_n$ in terms of the 
diagonal of a polygon. This was later adapted to type D$_n$ by 
Schiffler in~\cite{schi}, using a punctured $n$-gon. 
In the second part, we explain how to describe the $m$-cluster 
category of type A$_n$ in terms of so-called $m$-diagonals 
in a polygon. 

%%%%%%%%%%%%%%%%%%%%
\subsection{Quiver of diagonals}
%%%%%%%%%%%%%%%%%%%%

Let $\Pi$ be a polygon with $n+2$ vertices, 
labeled clockwise by $\{1,2,\dots,n+2\}$. To it, we 
associate a quiver $\Gamma(n,1)$ as follows: 

The vertices of $\Gamma(n,1)$ are the diagonals 
$\{(i,j)\mid 1\le i,j\le n+2,\ |i-j|>1\}$ of $\Pi$. 
The arrows are $(i,j)\to (i,j+1)$ and $(i,j)\to (i+1,j)$, 
where $i+1$ and $j+1$ are taken 
modulo $n+2$, 
provided the image 
is also a diagonal of $\Pi$. 
Furthermore, we define a bijection $\tau$ on the vertices of 
$\Gamma(n,1)$ 
as follows: $\tau:\,(i,j)\mapsto (i-1,j-1)$ (again, taking $i-1$ and $j-1$ 
modulo $n+2$). The quiver $\Gamma(n,1)$ together with this map 
$\tau$ is a stable translation quiver (cf. definition in 
Section~\ref{s:powers} below). 

\begin{example}[Hexagon]
Let us illustrate this in the case $n=4$. 
% 
%\begin{figure}
\begin{center}
  \psfragscanon
   \psfrag{1}{$1$}
   \psfrag{2}{$2$}
   \psfrag{3}{$3$}
   \psfrag{4}{$4$}
   \psfrag{5}{$5$}
   \psfrag{6}{$6$}
   \psfrag{(1,3)}{$(1,3)$}
   \psfrag{(1,4)}{$(1,4)$}
   \psfrag{(3,6)}{$(3,6)$}
\includegraphics[scale=.4]{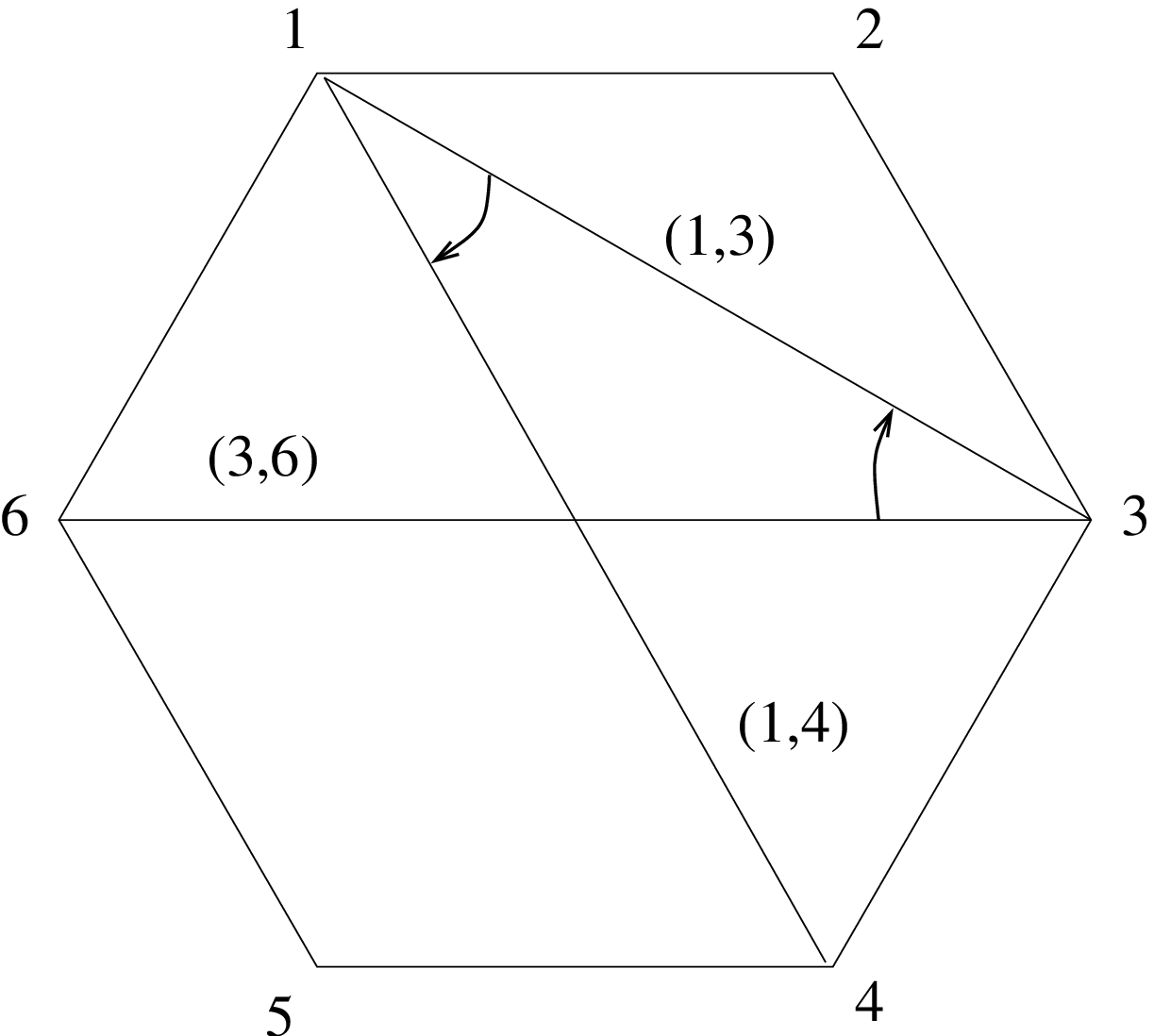}
\end{center}
The translation quiver $\Gamma(4,1)$ obtained from the hexagon is:  
$$
\xymatrix@-7mm{ 
 & & (1,5)\ar@{<.}[rr] \ar[rdd] & & (2,6) \ar[rdd]\ar@{<.}[rr] \ar[rdd] & & (1,3)\ar[rdd]\\
 \\
 & (1,4)\ar@{<.}[rr] \ar[rdd] \ar[ruu] & & (2,5)\ar@{<.}[rr]\ar[rdd]\ar[ruu] 
 & & (3,6)\ar[rdd] \ar@{<.}[rr] \ar[ruu] & & (1,4)\ar[rdd]\\
 \\
 (1,3)\ar@{<.}[rr] \ar[ruu] & & (2,4)\ar@{<.}[rr]\ar[ruu] & & (3,5)\ar@{<.}[rr]\ar[ruu] 
 & & (4,6)\ar@{<.}[rr] \ar[ruu] & & (1,5)
} 
$$
We have repeated the first slice 
$(1,3)\to (1,4)\to (1,5)$ 
at the end to indicate how the quiver wraps around. 
\end{example}

Observe that the quiver $\Gamma(4,1)$ is equal to the Auslander-Reiten 
quiver of the 
cluster category of type A$_3$ (Figure~\ref{fig:A3-AR}). 
More generally, one can show that the quiver of 
diagonals of an $n+2$-gon encodes the cluster category of 
type A$_{n-1}$: 

\begin{theorem}[\cite{ccs}, Section 2]\label{thm:cluster-diag} 
Let $Q$ be a quiver of Dynkin type A$_{n-1}$. Then the Auslander-Reiten 
quiver of $\CC_Q$ is isomorphic to the quiver 
$\Gamma(n,1)$ of diagonals in an $n+2$-gon. 
\end{theorem}

As a consequence of this, the cluster category $\CC_Q$ is equivalent 
to the additive category generated by the mesh category of the stable 
translation quiver $\Gamma(n,1)$ of diagonals of an $n+2$-gon. 
For details, we refer the reader to~\cite{ccs} and \cite{bm}.

\begin{remark}
The cluster category of type D$_n$ can be modelled in a similar way 
if we use a punctured $n$-gon. This has been done by Schiffler in~\cite{schi} using 
arcs in the punctured polygon. 
\end{remark}

\begin{remark}
Every maximal collection of non-crossing diagonals of a polygon $\Pi$ 
(punctured or not) 
is a triangulation of $\Pi$. All maximal collections have the same number 
of elements, this number is an invariant of $\Pi$. It is called the 
{\em rank of the polygon}. The rank of an $n+2$-gon is $n-1$, the rank 
of a punctured $n$-gon is $n$. 
This leads us back to cluster algebras - for the connection between 
(punctured) polygons of rank $n$ and cluster algebras of type A$_n$ (of 
type D$_n$ respectively) we refer the reader to~\cite{fst}. 
\end{remark}

%%%%%%%%%%%%%%
%
\subsection{Quiver of $m$-diagonals}\label{ss:m-diag}
%
%%%%%%%%%%%%%%

In a similar way, the $m$-cluster categories can be modelled using 
a certain %stable translation 
quiver $\Gamma(n,m)$. We will now explain 
how this works. 
Let $\Pi$ be an $nm+2$-gon. 
The vertices of $\Gamma(n,m)$ are the $m$-diagonals, i.e. the diagonals of 
the form $(1,m+2)$, $(1,2m+2)$, etc., where vertices are taken modulo 
$nm+2$. More precisely, an $m$-diagonal divides 
$\Pi$ into an $mj+2$-gon and an $m(n-j)+2$-gon (for $1\le j\le\frac{n-1}{2}$). 
The arrows send $(i,j)$ to $(i,j+m)$ and to $(i+m,j)$ whenever the image 
is also an $m$-diagonal. 
Furthermore, we define a translation $\tau_m$ on $\Gamma(n,m)$ which 
sends $(i,j)$ to $(i-m,j-m)$ (taking vertices modulo $nm+2$). 
Then $\Gamma(n,m)$ is also a stable translation quiver. 
With $m=1$ we just recover the case of the usual diagonals as 
described above. 

\begin{example}\label{ex:8gon}
Let $m=2$ and $n=3$; $\Pi$ is an octagon. 
Its $2$-diagonals 
are of the form $(1,4)$, $(1,6)$, $(2,5)$, etc.

%\begin{figure}
\begin{center}
  \psfragscanon
   \psfrag{1}{$1$}
   \psfrag{2}{$2$}
   \psfrag{3}{$3$}
   \psfrag{4}{$4$}
   \psfrag{5}{$5$}
   \psfrag{6}{$6$}
   \psfrag{7}{$7$}
   \psfrag{8}{$8$}
   \psfrag{(1,4)}{$(1,4)$}
   \psfrag{(1,6)}{$(1,6)$}
   \psfrag{(4,7)}{$(4,7)$}
\includegraphics[scale=.5]{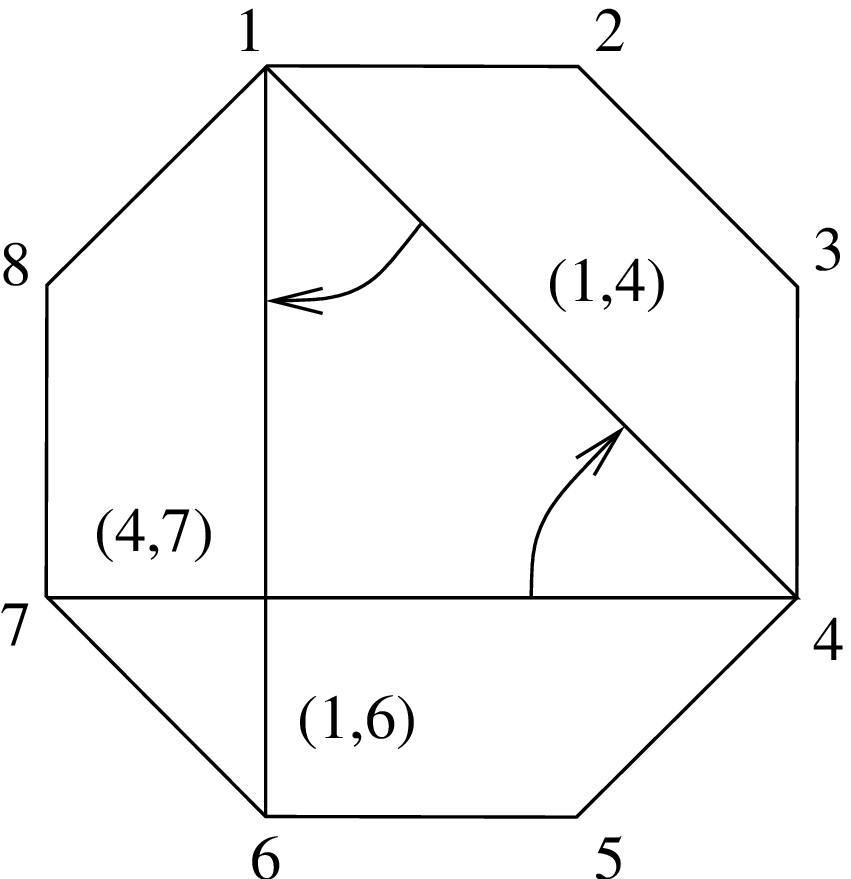}
\end{center}
Then the quiver $\Gamma(3,2)$ is: 
$$
\xymatrix@-8mm{ 
 & (1,6)\ar@{<.}[rr] \ar[rdd]  & & (3,8)\ar@{<.}[rr]\ar[rdd]
 & & (2,5)\ar[rdd]\ar@{<.}[rr]   
 & & (4,7)\ar@{<.}[rr] \ar[rdd] & & (1,6)  \\ 
  \\
 (1,4)\ar@{<.}[rr] \ar[ruu] & & (3,6)\ar@{<.}[rr]\ar[ruu] & & (5,8)\ar@{<.}[rr]\ar[ruu] 
 & & (2,7)\ar@{<.}[rr]\ar[ruu]   & & (1,4) \ar[ruu] 
} 
$$
Here, we have also repeated the first slice $(1,4)\to (3,6)$ to indicate 
how the quiver wraps around. 
\end{example}

Observe that the quiver $\Gamma(3,2)$ is just 
the Auslander-Reiten quiver of the 
$2$-cluster category of type A$_2$. So the $2$-diagonals in the octagon 
model the cluster category $\CC_{Q}^2$ for $Q$ of Dynkin type A$_2$. 
This holds more generally by the following result. 

\begin{theorem}[\cite{bm}]\label{thm:mcluster-mdiag} 
Let $Q$ be a quiver of Dynkin type A$_{n-1}$, let $m\ge 1$. 
Then the Auslander-Reiten 
quiver of $\CC_Q^m$ is isomorphic to the quiver 
$\Gamma(n,m)$ of $m$-diagonals in an $nm+2$-gon. 
\end{theorem}
Note that we recover Theorem~\ref{thm:cluster-diag} in the case $m=1$.

\begin{remark}
(1) Theorem~\ref{thm:mcluster-mdiag} tells us that in order to understand the 
$m$-cluster category of type A$_n$ it is enough to study $\Gamma(n,m)$. 

(2) 
To model the $m$-cluster categories of type D$_n$, one defines so-called 
$m$-arcs in a punctured $nm-m+1$-gon and obtains a quiver 
$\Gamma_{\odot}(n,m)$. 
One can show that this is the Auslander-Reiten 
quiver of $\CC_Q^m$ where $Q$ is of type D$_n$, cf. Theorem 3.6 in~\cite{bm2}. 
\end{remark}

The maximal collections of non-crossing $m$-diagonals in an $nm+2$-gon 
(resp. in a punctured $nm-m+1$-gon) 
correspond to {\em $m+2$-angulations} of the polygon. The number of elements 
in such a maximal collection is again an invariant of the polygon. It 
is equal to $n-1$ (resp. to $n$).

%%%%%%%%%%%%%%
%
\section{Powers of translation quivers}\label{s:powers}
%
%%%%%%%%%%%%%%

We now provide another way of obtaining $m$-cluster categories. 
It uses the concept of the power of a translation quiver which was 
introduced in~\cite{bm}. In order to explain this, let us give the precise 
definition of a translation quiver. 

\begin{definition}
A {\em translation quiver} is a pair $(\Gamma,\tau)$ where 
$\Gamma$ is a quiver, possibly with infinitely many vertices and arrows; 
$\tau$ is an injective map from a subset of the vertices 
of $\Gamma$ to the vertices of $\Gamma$, 
such that the following holds: 
the number of arrows going from a vertex $x$ to $y$ 
equals the number of arrows from $\tau y$ to $x$ 
for all vertices $x,y$ of $\Gamma$. 
The map $\tau$ is called the {\it translation of $(\Gamma,\tau)$}. 

If $\tau$ is defined on all vertices (and thus bijective) then $(\Gamma,\tau)$ 
is a {\em stable} translation quiver. 
\end{definition}
We remark that in all examples of stable translation quivers appearing 
in this article, 
the number of arrows between $2$ vertices is always at most $1$. 

We recall that a composition 
$x_0\to x_1\to\dots\to x_{m-1}\to x_m$ of $m$ arrows $x_i\to x_{i+1}$ 
(where the $x_i$ are vertices of $\Gamma$) is a {\em path of length $m$}. 
Such a path 
is said to be {\it sectional} if $\tau x_{i+1}\neq x_{i-1}$ for
$i=1,\dots,m-1$ (for which $\tau x_{i+1}$ is defined), cf.~\cite{rin}. 

\begin{definition}
Let $(\Gamma,\tau)$ be a translation quiver. The {\em $m$-th power $\Gamma^m$ 
of $\Gamma$} is the quiver whose vertices are the same as 
the vertices of $\Gamma$ and whose 
arrows are the sectional paths of length $m$ of $\Gamma$. 
\end{definition}
One can show that if $(\Gamma,\tau)$ is a stable translation quiver, 
then the pair $(\Gamma^m,\tau^m)$ is also a stable translation quiver 
(\cite{bm}, Section 6). 
Note however that $(\Gamma^m,\tau^m)$ is not connected in general, 
even if $(\Gamma,\tau)$ is so. The following example illustrates 
this.

\begin{example}\label{ex:2-power}
Let $\Gamma$ be the quiver $\Gamma(6,1)$ of diagonals in an octagon, let 
$m=2$. The quiver of the octagon has five rows, with first slice 
$(1,3)\to (1,4)\to (1,5)\to (1,6)\to (1,7)$:  

$$
\xymatrix@-9mm{ 
 &&& & (1,7)\ar@{<.}[rr] \ar[rdd]  & & (2,8)\ar[rdd]\ar@{<.}[rr]  && (1,3)\ar[rdd]
  \\ 
  \\
 && & (1,6)\ar@{<.}[rr] \ar[rdd]\ar[ruu]  & & (2,7)\ar@{<.}[rr]\ar[rdd]\ar[ruu]
 & & (3,8)\ar[rdd]\ar@{<.}[rr] \ar[ruu]  
 & & (1,4) \ar[rdd]   \\ 
  \\
 & & (1,5)\ar@{<.}[rr] \ar[rdd]\ar[ruu]  & & (2,6)\ar@{<.}[rr]\ar[rdd]\ar[ruu]
 & & (3,7)\ar[rdd]\ar@{<.}[rr]\ar[ruu]   
 & & (4,8)\ar@{<.}[rr] \ar[rdd]\ar[ruu] & & (1,5)\ar[rdd]  \\ 
  \\
 & (1,4)\ar@{<.}[rr] \ar[rdd]\ar[ruu]  & & (2,5)\ar@{<.}[rr]\ar[rdd]\ar[ruu]
 & & (3,6)\ar[rdd]\ar@{<.}[rr]\ar[ruu]   
 & & (4,7)\ar@{<.}[rr] \ar[ruu]\ar[rdd] & & (5,8)\ar@{<.}[rr]\ar[rdd]\ar[ruu]  
 & & (1,6)\ar[rdd]  \\ 
  \\
 (1,3)\ar@{<.}[rr] \ar[ruu] & & (2,4)\ar@{<.}[rr]\ar[ruu] & & (3,5)\ar@{<.}[rr]\ar[ruu] 
 & & (4,6)\ar@{<.}[rr]\ar[ruu]   & & (5,7)\ar@{<.}[rr]  \ar[ruu] 
 && (6,8)\ar[ruu]\ar@{<.}[rr]  && (1,7)
} 
$$
As before, we repeat the first slice at the end to indicate how the quiver 
wraps around. 

The second power of $\Gamma(6,1)$ has three components. One containing 
the vertex $(1,3)$, one containing the vertex $(1,4)$ and the third containing 
the vertex $(2,4)$: 

$$
\xymatrix@-8mm{ 
 & (1,6)\ar@{<.}[rr] \ar[rdd]  & & (3,8)\ar@{<.}[rr]\ar[rdd]
 & & (2,5)\ar[rdd]\ar@{<.}[rr]   
 & & (4,7)\ar@{<.}[rr] \ar[rdd] & & (1,6)  \\ 
  \\
 (1,4)\ar@{<.}[rr] \ar[ruu] & & (3,6)\ar@{<.}[rr]\ar[ruu] & & (5,8)\ar@{<.}[rr]\ar[ruu] 
 & & (2,7)\ar@{<.}[rr]\ar[ruu]   & & (1,4) \ar[ruu] 
} 
$$

$$
\xymatrix@-8mm{ 
 && (1,7)\ar@{<.}[rr] \ar[rdd]  & & (1,3)\ar[rdd]
 \\ 
  \\
 & (1,5)\ar@{<.}[rr] \ar[ruu]\ar[rdd] & & (3,7)\ar@{<.}[rr]\ar[ruu]\ar[rdd] 
  & & (1,5)\ar[rdd] 
 & &  \\
 \\
 (1,3)\ar@{<.}[rr]\ar[ruu]  & & (3,5)\ar@{<.}[rr]\ar[ruu] 
  & & (5,7)\ar@{<.}[rr]\ar[ruu]  & & (1,7) 
} 
$$
and 
$$
\xymatrix@-8mm{ 
 && (2,8)\ar@{<.}[rr] \ar[rdd]  & & (2,4)\ar[rdd]
 \\ 
  \\
 & (2,6)\ar@{<.}[rr] \ar[ruu]\ar[rdd] & & (4,8)\ar@{<.}[rr]\ar[ruu]\ar[rdd] 
  & & (2,6)\ar[rdd] 
 & &  \\
 \\
 (2,4)\ar@{<.}[rr]\ar[ruu]  & & (4,6)\ar@{<.}[rr]\ar[ruu] 
  & & (6,8)\ar@{<.}[rr]\ar[ruu]  & & (2,8) 
} 
$$
\end{example}

We observe that the component through the vertex $(1,4)$ is the same as the 
quiver $\Gamma(3,2)$ from Example~\ref{ex:8gon}. This is a property 
that holds in general, cf. Theorem~\ref{thm:conn-comp}. 
The other two components are isomorphic to the Auslander-Reiten quiver 
of the orbit category $\DD^b(kQ)/[1]$ where $Q$ is of Dynkin type A$_3$.

\begin{theorem}\label{thm:conn-comp}
The quiver $\Gamma(n,m)$ is a connected component of $\Gamma(nm,1)^m$. 
\end{theorem}

From Theorem~\ref{thm:conn-comp} one obtains that 
the $m$-cluster category of type A$_n$ is a full subcategory 
of the additive category generated by the mesh category of 
$\Gamma(nm,1)$ (for a definition of the mesh category we refer the reader 
to~\cite[Section 3]{bm}). 
In other words: in order to understand the $m$-cluster category, 
there is an alternative approach to the one presented in Subsection~\ref{ss:m-diag}. 
Namely, we can consider the $m$-th power 
of the quiver given by the usual diagonals. The Auslander-Reiten quiver 
of the $m$-cluster category of type A$_{n-1}$ appears as a connected 
component of $\Gamma(nm,1)^m$. 

\begin{remark}
We can actually prove an analogous result as Theorem 
~\ref{thm:conn-comp} for type D$_n$, cf.~\cite{bm2}: 
Let $\Gamma_{\odot}(n,1)$ denote the quiver obtained from the arcs of a punctured 
$n$-gon (note the difference: here, the polygon has $n$ vertices instead of 
$n+2$) and by 
$\Gamma_{\odot}(n,m)$ the quiver of $m$-arcs in a punctured $nm-m+1$-gon. 
Then the Auslander-Reiten quiver of the $m$-cluster category of type 
D$_n$ is $\Gamma_{\odot}(n,m)$ and it is a connected component of 
$\Gamma_{\odot}(nm-m+1,1)^m$. 
For details we refer to Section 5 of~\cite{bm2}\footnote{to be precise: 
the arrows in the $m$-power arise from {\em restricted} sectional paths of length $m$, 
and we are taking the restricted $m$th power.}. 
\end{remark}

A natural question to ask at this point is what the other connected components 
of the $m$th power of $\Gamma(nm,1)$ are. 
Remarkably, this question is much easier in type D. There, we have 
a complete answer: 

The connected components of the (restricted) $m$th power of 
the Auslander-Reiten quiver of the cluster category of type D$_{nm-m+1}$ 
are exactly the union of the Auslander-Reiten quiver of the $m$-cluster 
category of type D$_n$ with $m-1$ copies of the Auslander-Reiten 
quiver of $\DD^b(\mbox{A}_{n-1})/\tau^{nm-m+1}$, cf. 
Theorem 5.2 of~\cite{bm2}. 

The difficulty arising in type A has to do with the additional symmetry of the 
Dynkin graph of type A, i.e. with the involution sending the first to the 
last node, the second to the second-to-the-last, etc. 

Let now $Q$ be of Dynkin type A$_n$. For odd $m$, 
C. Ducrest(\cite{du08}) was recently able to answer the question of the components 
of the $m$th power. For even $m$, she provides a partial description. 
Her result is the following (see Section 4 of~\cite{du08}): 

\begin{theorem} 
$$
\Gamma(nm,1)^m=\Gamma(n,m)\ \ \bigsqcup\ \ \bigcup_{i=1}^t \Gamma_i
$$ 
where $\Gamma_i$ is the Auslander-Reiten quiver of an orbit category 
of $\DD^b(\mbox{A}_n)$ under an auto-equivalence of the form 
$\tau^{-s}\circ[r]$ for some $t$ and for some $s$, $r$ where we 
can assume $s<n$. 

Furthermore, if $m$ is odd then $r=(m-1)/2$ and $s=\frac{1}{2}((m-1)(n-1))+1$. 
\end{theorem}

The even case is more complicated. One can show that 
for even $m$, we have $m/2\le r\le m$. Example~\ref{ex:2-power} above shows 
that $r=m$ does occur. 

\begin{remark}
If $m$ is odd, one can show that the $m$th power only has one connected 
component per row of the original quiver $\Gamma(nm,1)$. 
In the even case, there are examples where we obtain one component per 
row and examples where there are two components per row. 
For details we refer to~\cite{du08}. 
\end{remark}

C. Ducrest has developed a programme to calculate all components 
of the $m$th power of $\Gamma(m,n)$ for all $n,m\le 20$. This programme 
is available online at \\
http://www.math.ethz.ch/$\sim$baur/algo/ \\
and the documentation explaining how it works is ~\cite{du-doc}. 
We hope that this programme 
will help us finding the complete answer for $m$ even. 

%
%
%%%%%%%%%%%%%%%%%%%%%%%%%%%%

\end{document}